\documentclass[reqno,12pt]{amsart}
\usepackage[utf8]{inputenc}
\usepackage{indentfirst}
\usepackage{amsmath}
\usepackage[T1]{fontenc}
\usepackage{amssymb}			
\usepackage{amsthm}
\usepackage{amsfonts}
\usepackage{enumerate}
\usepackage{enumitem}
\usepackage{setspace}
\usepackage{mathtools}
\usepackage[body={7in,9.5in},top=0.5in, left=0.8in ]{geometry}
\usepackage{color}
\newcommand{\bn}{\mathbb{N}}

\newcommand{\ov}[1]{\overline{#1}}

\newcommand{\pref}[1]{\textup{(\ref{#1})}}

\newtheorem{thrm}{Theorem}%[section] %definicje wlasne; kolejno

\newtheorem{thm}[thrm]{Theorem}
\newtheorem{definition}[thrm]{Definition}
\newtheorem{remark}[thrm]{Remark}

\newtheorem{lm}[thrm]{Lemma}

\newtheorem{rmk}[thrm]{Remark}

\newtheorem{xmpl}[thrm]{Example}
\newtheorem{theorem}[thrm]{Theorem}

\usepackage[nobysame,initials]{amsrefs}

\newcommand{\X}{{\mathbb{X}}}

\newcommand{\K}{{\mathbb{K}}}
\newcommand{\R}{{\mathbb{R}}}
\newcommand{\C}{{\mathbb{C}}}
\newcommand{\Z}{{\mathbb{Z}}}
\newcommand{\N}{{\mathbb{N}}}

\title[On composition and Right Distributive Law for formal power series]{On composition and Right Distributive Law for formal power series of multiple variables}

\author{Dariusz Bugajewski}
\author{Alessia Galimberti}
\author{Piotr Ma\'ckowiak}

\address[D. Bugajewski, P. Ma\'ckowiak]{Department of Nonlinear Analysis and Applied Topology\\
  Faculty of Mathematics and Computer Science\\
  Adam Mickiewicz University\\
  Uniwersytetu Pozna\'nskiego 4\\
  61-614 Pozna\'n\\
  Poland}

\address[A. Galimberti]{Dipartimento di Matematica\\
Universit\`a degli Studi di Milano\\
Via Cesare Saldini~50\\20133 Milano\\
Italy}

\email[D.~Bugajewski]{ddbb@amu.edu.pl}
\email[A.~Galimberti]{alessia.galimberti1@studenti.unimi.it}
\email[P.~Ma\'ckowiak]{piotr.mackowiak@amu.edu.pl}

\begin{document}

\begin{abstract} In the first part of the paper we prove a necessary and sufficient condition for the existence of the composition of formal power series in the case when the outer series is a series of one variable while the inner one is a series of multiple variables. The aim of the second part is to remove ambiguities connected with the Right Distributive Law for formal power series of one variable as well as to provide analogues of that law in the multivariable case. 
\end{abstract}

\maketitle

\section{Introduction}

There is no need to convince about how important {is the role that formal power series, as well as formal Laurent series, play} in various branches of mathematics. In particular, they are frequently used in solving various type {of} equations, like ordinary differential equations or partial differential equations (see e.g. \cite{PraSp}, \cite{Sib} among numerous papers concerning this topic). Other applications of formal power series to Riordan groups and proofs of some classical results, like the Cayley-Hamilton Theorem, can be found e.g. in \cite{CMP}, \cite{Gan M2} and \cite{JLN}. 

An elementary introduction to the theory of formal power series in several variables can be found e.g. in the paper \cite{Hau}. Moreover, in that paper one can find further applications of formal power series to the theory of multiplicative arithmetical functions of several variables and to the theory of the cumulants of random vectors. 

In this article we are going to discuss two basic topics appearing in the theory of formal power series, that is, the composition of formal power series and the Right Distributive Law.

Our first goal is to prove a necessary and sufficient condition for the existence of the composition of formal power series in the case when {the} outer series is a series of one variable while {the} inner one is a series of multiple variables. It appears that the condition of our theorem is similar to the one appearing in Gan-Knox's Theorem involving the composition of two formal power series of one variable (see \cite {GK} for more details). 

Our second goal is to examine the General Right Distributive Law. First, we are going to explain ambiguities connected with this law in the one variable case. Second, we prove analogues of that law for formal power series of multiple variables.

{The reader can find more information concerning the composition of formal power series}, e.g., in \cite{BG} or \cite{Daw}, while more details concerning the Right Distributive Law can be found, e.g., in \cite{Gan1}.

\section{Definitions, conventions and some basic facts}

The goal of this section is to recall and introduce basic definitions and facts which will be needed in the sequel. 

The symbol $\Z$ denotes the set of integer numbers, $\N$ is the set of positive integers and $\N_0:=\N\cup\{0\}$. If $n\in \N$, then $[n]:=\{1,\ldots,n\}$. Let us fix $q\in \N$, $k\in \N_0$ and denote by $C_k$ (usually the value of $q$ will be fixed and clear from the context, so we write just $C_k$ instead of, say, $C^q_k$)  the set of all nonnegative integer solutions $c_1,\ldots,c_q$ of the equation $c_1+\ldots+c_q=k$ for $k\in \N_0,\, q\in \N$, that is 
$$
C_k:=\{c=(c_1,\ldots, c_q)\in \N_0^q:\, c_1+\ldots+c_q=k\}.
$$ 
It is known (see \cite{S}, p.25) that the number of elements of $C_k$ is $\binom{q+k-1}{q-1}$ - let us denote that number by $s_k$. It is clear that $s_0=1$ and $C_k\cap C_{k'}=\emptyset $ for $k\neq k'$. 
Obviously, $C:=\N^q_0=\bigcup_{k\in \bn_0}C_k$ and  for each $c\in C$ there is exactly one $k\in \bn_0$ for which $c\in C_k$; we denote such a $k$ by $k(c)$.
Let $\K$ stand for the field of real (or complex) numbers. 

Now, we are ready to state the definition of a formal power series of multiple variables.

\begin{definition}
A formal power series $f$ in $q$-variables $x:=(x_1,\ldots,x_q)\in \K^q$ (for short: q-fps) is a formal sum of the form
$f(x):=\sum_{c\in C}f_cX^c$, where $f_c\in \K$ and $X^c:=x_1^{c_1}\ldots x_q^{c_q}$ for all $c\in C$. An element $f_c$ is called the $c$-th coefficient of the q-fps $f$, $c\in C$. The set of all q-fps is denoted by $\X_q(\K)$ or just by $\X_q$.
\end{definition}

\begin{rmk} \textup{
Let us observe that a q-fps $f$ can be uniquely identified with the mapping $C\ni c\mapsto f_c\in\K$.} 
\end{rmk}

\begin{rmk}\textup{ For a q-fps with $q=1$ we equate $x=(x_1)$ with $x$, where $x$ is the name of the variable, so we write $x^c$ instead of $X^c$. Further, for $f\in\X_1$, we identify $\sum_{c\in \N_0}f_cx^c$ with $\sum_{n=0}^{\infty}f_nx^n$ as long as we do not assign any particular value to $x$, that is, as long as we do not consider the value of the sum $\sum_{c\in \N_0}f_cx^c$ corresponding to a specific value of $x$.}
\end{rmk}

We distinguish two special q-fps, denoted by $I$ and $\theta$, whose coefficients are defined as $I_{(0,\ldots,0)}:=1$, and $I_c:=0$ for $c\in C\setminus\{(0,\ldots,0)\}$, where $(0,\ldots,0)\in C$, $0$ is the additive identity of the field $\K$, $1$ is the multiplicative identity of the field $\K$, and $\theta_c:=0$, $c\in C$, respectively. Let us observe that the q-fps $I$ and $\theta$ are neutral elements of the multiplication and the addition of q-fps, respectively.

For q-fps $f,\,g$, the sum $f+g$, difference $f-g$ and product $af$, where $a\in \K$, are q-fps whose coefficients are defined as $f_c+g_c, \, f_c-g_c$, and $af_c$, $c\in C$, respectively. We write $f=g$ if and only if $f_c=g_c$ for every $c\in C$. Obviously $f=g$ if and only if $f-g=\theta$. 

Let us recall the definition of the Cauchy product of q-fps. 

\begin{definition}
For q-fps $f(x)=\sum_{c\in C}f_cX^c, \, g(x)=\sum_{c\in C}g_cX^c\in \X_q$, $q\in \bn$, the $q$-dimensional Cauchy product of $f$ and $g$ is a q-fps $h=fg$ defined as
$$h(x):=\sum_{c\in C}\underbrace{\left(\sum_{a,\,b\in C:\,a+b=c}f_ag_b\right)}_{h_c:=}X^c.$$
\end{definition}

Obviously, there are a finite number of pairs $(a,b)$ such that $a, b\in C$ and $a+b=c$ for the given $c\in C$. It is also clear that for such a pair we have $a\leq c$ and $b\leq c$, where the inequality $\leq$ is understood coordinatewise--like. %Moreover, if $a,b,c\in C$ and $a+b=c$, then $k(a)\leq k(c)$ and $k(b)\leq k(c)$.

Now, let us introduce some further conventions we shall use. For any $a\in \K$ we put $a^0:=1$ even if $a=0$. For a fixed $c\in C$, the monomial $X^c=x_1^{c_1}\ldots x_n^{c_q}$, where $x=(x_1,\ldots,x_q)\in \K^q$, can be treated as a q-fps with all, but $c$-th, coefficients equal to $0$. Moreover, we equate $X^c$ with $1X^c$. Any finite formal sum $a_1X^{c_1}+\ldots+a_nX^{c_n}$, $c_1,\ldots, c_n\in C$, $a_1,\ldots,a_n\in \K$, can be viewed as a q-fps whose all coefficients, possibly except $c_i$-th, $i\in [n]$, are $0$, and $c_i$-th coefficient is $a_i$, $i\in[n]$. This way any element $a\in \K$ can be identified with the constant q-fps $aX^{(0,\ldots,0)}=aI(x)$. 
\textcolor{green}{} It is also obvious that the  addition of q-fps is associative and commutative. 
Basic properties of the formal product are given in the following lemma. Its proof is standard so we omit it. 

\begin{lm}\label{lm:0}
Let $f,g,h\in \X_q$ be fixed for a given $q\in \bn$. Then: 
\begin{enumerate}[label={\textup{(\arabic*)}},ref=\textup{(\arabic*})]
\item\label{thrm:0:1} $fg=gf$,
\item\label{thrm:0:2} $fI=If=f$,
\item\label{thrm:0:3} $(fg)h= f(gh)$,
\item\label{thrm:0:4} $(f+g)h= (fh)+ (gh)$,
\item\label{thrm:0:5} $f\theta=\theta$,
\item\label{thrm:0:6} $fg=\theta$, implies $f=\theta$ or $g=\theta$ (which means that $\X_q$ is an integral domain).
\end{enumerate}
\end{lm}

We are now ready to define the composition of a $1$-fps $g$ with a $q$-fps $f$. 

\begin{definition}\label{df:composition}
For $f(x)=\sum_{c\in C}f_cX^c\in \X_q$, $q\in \bn$, and  $g(y)=\sum_{n=0}^{\infty}g_ny^n \in\X_1$, the composition of $g$ with $f$ is a $q$-fps $h=g\circ f$ defined by
$$h(x):=\sum_{c\in C}\underbrace{\left(\sum_{n=0}^\infty g_n f^{n}_c\right)}_{h_c:=}X^c,$$
provided that the coefficients $h_c$ exist, that is, if the series defining $h_c$ converge for every $c\in C$, where $f^{n}\in \X_q$ is the $n$-th power of $f$, that is, $f^{n}:=\underbrace{ff\ldots f}_{n\times}$, $f^n_c:=(f^n)_c$, $n\in \N$, $c\in C$, and $f^{0}:=I$.
\end{definition}

\begin{remark}\label{order}\textup{
Given a q-fps, if one needs an order of summation of its terms, then one can proceed in the following way.  
Let us define a linear order $\preceq $ on $C$ as follows: for $c=(c_1,\ldots, c_q),\ c'=(c'_1,\ldots, c'_q)\in C$  $$c'\preceq c\Leftrightarrow (c'_1+\ldots+c'_q<c_1+\ldots+c_q \vee c\geq_{lex}c'),$$
where $c\geq_{lex}c'$ means that $c$ is lexicographically not less than $c'$. For each $k\in \N$, let us denote by $(c^k_1,\ldots,c^k_{s_k})$ the sequence of all elements of $C_k$ satisfying $c^k_{i+1}\geq_{lex}c^k_{i},\, i \in [s_k -1]$.
Define
$$
\ov{C}:=(c^0_1,c^1_1,\ldots,c^1_{s_1}, c^2_1,\ldots,c^2_{s_2},\ldots),
$$
 where $c^k_i\in C$, $c^k_i\preceq c^{k'}_{i'}$ whenever $k<k'$ or $k=k' \wedge i< i'$, and $c^k_i=c^{k'}_{i'}$ if and only if $i=i'\, \wedge\, k=k'$. Since all the terms of $\ov{C}$ are different, one can identify $C$ with $\ov{C}$ and we can write $C=\{c^0_1,c^1_1,\ldots,c^1_{s_1}, c^2_1,\ldots,c^2_{s_2},\ldots\}$. The sequence $\ov{C}$ induces the order of addition in the sum $\sum_{c\in C}f_cX^c$ in a natural way.}
\end{remark}

\begin{definition}
For a q-fps $f\in \X_q$ and $k\in \bn_0$, the $k$-th block of $f$ is a q-fps $f[k]\in \X_q$ whose coefficients are given by $f[k]_c:=f_c$, $c\in C_k,$ and $f[k]_c:=0$, otherwise.
\end{definition}

As the sets $C_k$, $k\in \bn_0$, are finite and pairwise disjoint, we get $f_c=\sum_{k\in \N_0}f[k]_c=\sum_{k=0}^{\infty}f[k]_c=f[k(c)]_c$ for $c\in C$. Further, $f(x)=\sum_{c\in C}\left(\sum_{k=0}^\infty f[k]_c\right)X^c=\sum_{k=0}^\infty \left(\sum_{c\in C}f[k]_cX^c\right)=\sum_{k=0}^\infty f[k](x)$, under the convention that the (formal) sum of countably many terms of the form $0X^c$, $c\in C$, is $\theta$.

For the coefficients of the $n$-th power of a $q$-formal power series $f$, we have an obvious:

\begin{lm}\label{lm:powers}
For a q-fps $f(x)=\sum_{c\in C}f_cX^c\in \X_q$, $q\in \bn$, and any $n\in \N$ and $c\in C$, it holds $f^n_c=\sum f_{c_1}\ldots f_{c_n}$, where the sum is taken over all solutions of the equation $c_1+\ldots+c_n=c$, $c_i\in C,\,i\in [n]$.
\end{lm}
\begin{proof} Let $c\in C$, $n\in \bn$ be fixed. If $n=1$, then the claim is true. If $n>1$, by associativity of multiplication we have $f^n_c=\sum_{a+c_n=c}f^{n-1}_af_{c_n}=\sum_{a+c_n=c}(\sum_{b+c_{n-1}}f^{n-2}_bf_{c_{n-1}})f_{c_n}=\sum_{a+c_{n-1}+c_n=c}f^{n-2}_af_{c_{n-1}}f_{c_n}=\ldots=\sum_{c_{1}+\ldots+c_n=c}f_{c_1}\ldots f_{c_n}$, where all elements in equalities under summation symbols are members of $C$.
\end{proof} 

By Lemma \ref{lm:powers}, we conclude that for given $c\in C$ and $n\in \N$, we have $$f^n_c=\sum_{c_{1}+\ldots+c_n=c}(f[0]+\ldots+f[k])_{c_1}\ldots (f[0]+\ldots+f[k])_{c_n},$$ where $k:=k(c)$. This implies that to compute the $c$-th coefficient of $f^n$ it is enough to compute the $c$-th coefficient of the q-fps $(f[0]+\ldots+f[k])^n$, where $k:=k(c)$. Since formal addition and multiplication are associative and commutative, for a  given $c\in C$, by the Multinomial Theorem, we get $(f[0]+\ldots +f[k])^n=\underbrace{(f[0]+ \ldots + f[k])\ldots (f[0]+ \ldots + f[k])}_{n\times}=\sum \frac{n!}{v_0! \ldots  v_k!}f[0]^{v_0}\ldots f[k]^{v_k}$, where $k:=k(c)$ and the sum extends over nonnegative integer solutions of $v_0+v_1+\ldots+v_k=n$. Observe that the only components of the sum that contribute to the $k$-th block of $f^n$, that is, $f^n[k]$, are those for which $v_1+2v_2+\ldots+kv_k=k$. Therefore, for a q-fps $f$ and any $c\in C$, $n\in \N$, it follows that \begin{equation}\label{eqn:powers} f^n_c=
\sum \frac{n!}{v_0! \ldots  v_k!}(f[0]^{v_0}\ldots f[k]^{v_k})_c,\end{equation} where $k=k(c)$ and the sum extends over nonnegative integer solutions of the system \begin{equation}\label{eqn:cndtns}\left\{\begin{array}{l}v_0+v_1+\ldots+v_k=n\\v_1+2v_2+\ldots+kv_k=k.\end{array}\right.\end{equation}

\section{Existence of composition of two formal power series}

In this section, we provide {a} proof of a necessary and sufficient condition for the existence of composition of $g\in \X_1$ with $f\in \X_q,\, q\in \N$. We use the following conventions:  $\sum_{i\in \emptyset}:=0$, $\sum_{i\in A}0:=0$ ($A$ is a countable set) and $b^0:=1$ for $b\in \K$.

\begin{thm}\label{thm:gencomp} Fix $q\in \N$. Let $g(y):= \sum_{n\in\N_0}^\infty g_ny^n\in \X_1$, $f(x):=\sum_{c\in C}f_cX^c\in \X_q$ and $f\neq bI$ for any $b\in \K$. The composition $h:=g\circ f$ exists if and only if \begin{equation}\label{eqn:conv}
    \sum_{n=k}^\infty \binom{n}{k} g_nb_0^{n-k}\in \K,
  \end{equation}
where $b_0=f_{(0,\ldots,0)}$, for $k\in \N_0$.
\end{thm}
\begin{rmk}\textup{
Before we present a proof of Theorem \ref{thm:gencomp}, let us notice that {if we omit the assumption that $f\neq bI$, then $g\circ f$ exists if and only if the series $\sum_{n=0}^\infty g_nb^n$ converges}. In this case the composition is of the form $(g\circ f)(x)=(\sum_{n=0}^\infty g_nb^n)X^{(0,\ldots,0)}$.}
\end{rmk}
\begin{proof}
The main line of this proof is similar to the proof of Theorem 3 {in} \cite{BM}. By formulas (\ref{eqn:powers})-(\ref{eqn:cndtns}) and Definition \ref{df:composition}, it follows that the composition $h$ of $g$ with $f$ exists if and only if, for every $k\in \N_0$, $c\in C_k$ there exists the limit
\begin{equation}\label{eqn:lim}
\lim_{d\to \infty}\sum_{n=0}^{d}g_n(w[k]^n)_c\in \K,\end{equation} 
where $w[k]^n:=\sum \frac{n!}{v_0! \ldots  v_k!}f[0]^{v_0}\ldots f[k]^{v_k}$ and the sum is taken over the set of nonnegative integer solutions $v=(v_0,v_1,\ldots,v_k)$ to the system (\ref{eqn:cndtns}).
 For each $k\in \N_0$, the existence of the limit \pref{eqn:lim} is equivalent to the existence of the limit \begin{equation}\label{eqn:lim1}\lim_{d\to \infty}\sum_{n=k}^{d}g_n(w[k]^n)_c\in \K.
\end{equation}
Let us denote $w_k^{d}(x):=\sum_{n=k}^{d}g_nw[k]^n(x)$.
If $k=0$, then $(w^{d}_0)_{c_0}=\sum_{n=0}^{d}g_nb_0^n$, where $c_0:=(0,\ldots,0)\in C_0$ and the existence of $h_{c_0}$ is equivalent to the existence of the limit $\lim_{d\to \infty} \sum_{n=0}^{d}g_nb_0^n$.

Let $k\in \N_0$ and define 
$$R^n_k:=\left\{(v_0,\ldots,v_k)\in \N^{k+1}_0:\, \sum_{i=0}^kv_i=n\text{ and }\sum_{i=1}^kiv_i=k\right\}=\{v^{n,1},\ldots,v^{n,r^n_k}\},
$$
where $v^{n,j}=(v^{n,j}_0,v^{n,j}_1,\ldots,v^{n,j}_k)\in \N^{k+1}_0$ and $r^n_k$ denotes the number of elements of $R^n_k$,
$$m_{k,j}^n:=\frac{1}{v^{n,j}_1! \ldots v^{n,j}_k!}f[1]^{v^{n,j}_1}\ldots f[k]^{v^{n,j}_k},\quad j\in [r^n_k],$$
$$W^n_{k,s}:=\{j\in [r^n_k]:\, \sum_{i=1}^{k}v^{n,j}_i=s\},\quad s\in [k],$$
$$d^{n}_{k,s}:=\sum_{j\in W^n_s}m^n_{k,j}.$$
For $d\geq k$, it holds 
\begin{multline*}
    w^{d}_k=\sum_{n=k}^{d}g_n\left(\sum_{v\in R^n_k} \frac{n!}{v_0! \ldots  v_k!}f[0]^{v_0}\ldots f[k]^{v_k}\right)= \\ =\sum_{n=k}^{d}g_n\left(\sum_{j=1}^{r^n_k}b_0^{v^{n,j}_0}\frac{n!}{v^{n,j}_0!}\frac{1}{v^{n,j}_1! \ldots  v^{n,j}_k!}f[1]^{v^{n,j}_1}\ldots f[k]^{v^{n,j}_k}\right)=\\\sum_{n=k}^{d}g_n\left(\sum_{j=1}^{r^n_k}b_0^{v^{n,j}_0}\frac{n!}{v^{n,j}_0!}m^{n}_{k,j}\right)=\sum_{n=k}^{d}g_n\left(\sum_{s=1}^{k}\left(\sum_{j\in W^n_{k,s}}b_0^{v^{n,j}_0}\frac{n!}{v^{n,j}_0!}m^{n}_{k,j}\right)\right)=(\star)
\end{multline*}
but since for $n\geq k$ we have $W^n_{k,s}=W^{n'}_{k,s}$, $m^n_{k,j}=m^{n'}_{k,j}$, $d^{n}_{k,s}=d^{n'}_{k,s}$ as $n'\geq n$, there is no necessity to discern between $n$ and $n'$ as long as $n'\geq n\geq k$. Therefore we write $W_{k,s}$, $m_{k,j}$, $d_{k,s}$ instead of $W^{n}_{k,s}$, $m^n_{k,j}$, $d^n_{k,s}$, respectively. Due to the equality $v^{n,j}_0=n-v^{n,j}_1-\ldots-v^{n,j}_k$, we get
  \begin{multline}\label{eq:lim2}
  (\star)=\sum_{n=k}^{d}g_n\left(\sum_{s=1}^{k}\left(\sum_{j\in W_{k,s}}b_0^{v^{n,j}_0}\frac{n!}{v^{n,j}_0!}m_{k,j}\right)\right)=\\
  =\sum_{n=k}^{d}g_n\left(\sum_{s=1}^{k}\frac{n!}{(n-s)!}b_0^{n-s} d_{k,s}\right)=\sum_{s=1}^{k} s!\left(\sum_{n=k}^{d}{n \choose s} g_nb_0^{n-s}\right )d_{k,s}.
  \end{multline}
  If $b_0=0$, then the last formula reduces to $w^{d}_k=k!g_kd_{k,k}$ and the claim follows.
  It is obvious that validity of condition (\ref{eqn:conv}), for $k\in \N_0$, implies that $h$ is the composition of $g$ with $f$. This proves that condition \pref{eqn:conv} is sufficient for the existence of $h$.
	
To prove necessity part we assume that $h=g\circ f$ exists, so $\lim_{d\to \infty} (w^{d}_k)_c\in \K$, $c\in C$. Assume that $b_0\neq 0$. We shall proceed by induction with respect to $k$. As we know, the coefficient $h_{(0,\ldots,0)}$ is well--defined if and only if the series in (\ref{eqn:conv}) converges for $k=0$. Suppose that the equivalence is true for some $k\in \N_0$, which implies the convergence of series $\sum_{n=s}^\infty \binom{n}{s} g_nb_0^{n-s}$ in $\K$, $s=0,\,\ldots,k$. Let $l$ denote the least positive integer for which there is $c\in C_l$ with $f[l]_c\neq 0$ (in other words, $l$ is the least positive integer for which $f[l]\neq \theta$ - such a number exists since $f\neq b_0I$ and $b_0\neq 0$). Let us fix some positive integer $m$. First we will show that $d_{ml,s}=\theta$, whenever $ml\geq s>m$, and $d_{ml,m}=\frac{1}{m!}f[l]^m$ (recall that we consider $n,n'\geq k$, so $d^{n'}_{k,s}=d^{n'}_{k,s}=:d_{k,s}$). Indeed, assume that $v_1,\dots,v_{ml}\in \N_0$ solve~$v_1+2v_2+\dots+mlv_{ml}=ml$ and that $v_1+\dots+v_{ml}=s$. If $ml\geq s>m$, it cannot be true that $v_1=\dots=v_{l-1}=0$, since then we would have $ml=lv_l+\dots+mlv_{ml}\ge l(v_l+\dots+v_{ml})=sl>ml${,} a contradiction. Therefore, if $ml\geq s>m$, then in any component of the sum defining $d_{ml,s}$, there is some $i\in[l-1]$ such that $v_i>0$ and{,} due to the definition of $l${,} $f[i]^{v_i}=\theta$, so $d_{ml,s}=\theta$. Consider now the case when $s=m$. If for some $i\in [l-1]$ we have $v_i>0$, the component containing $f[i]^{v_i}$($=\theta$) adds nothing to $d_{ml,m}$, so we can restrict ourselves to the case where $v_1=\dots=v_{l-1}=0$. If there were any $i>l$ such that $v_i>0$, we would have $v_l+\dots+v_{ml}=m$ and $lv_l+\dots+mlv_{ml}=ml$; but $ml=l(v_l+\dots+v_{ml})<lv_l+\dots+mlv_{ml}=ml${,} again a contradiction. This shows that the only one component of in the sum defining $d_{ml,m}$ which is different from $\theta$ is $\frac{1}{m!}f[l]^m$, for any fixed $m$. In particular, 
$$d_{ml,m}(x)=\frac{1}{m!}f[l]^m(x).$$ 
For $m=k+1$ and $d\geq (k+1)l$, by \pref{eq:lim2}, we get
 \begin{multline*}
 	w^{d}_{(k+1)l}(x)=\sum_{s=1}^{(k+1)l} s!\left(\sum_{n=(k+1)l}^{d}{n \choose s} g_nb_0^{n-s}\right)d_{(k+1)l,s}(x)=\\ \sum_{s=1}^{k+1} s!\left(\sum_{n=(k+1)l}^d\binom{n}{s}g_nb_0^{n-s}\right)d_{(k+1)l,s}(x)=\\
	\left[\sum_{s=1}^{k} s!\left(\sum_{n=(k+1)l}^d\binom{n}{s}g_nb_0^{n-s}\right)d_{(k+1)l,s}(x)\right]+(k+1)!\left(\sum_{n=(k+1)l}^d\binom{n}{k+1}g_nb_0^{n-(k+1)}\right)d_{(k+1)l,(k+1)}(x)=\\
	\left[\sum_{s=1}^{k} s!\left(\sum_{n=(k+1)l}^d\binom{n}{s}g_nb_0^{n-s}\right)d_{(k+1)l,s}(x)\right]+\left(\sum_{n=(k+1)l}^d\binom{n}{k+1}g_nb_0^{n-(k+1)}\right)f[l]^{k+1}(x).
\end{multline*}
Let now $\ov{c}\in C_l$ be the unique minimal element of $\{c\in C_l:\, f[l]_c\neq0\}$ in the sense of the lexicographic ordering of $C_l\subset \N_0^q$. Then we have $(f[l]^{k+1})_{(k+1)\ov{c}}=f_{\ov{c}}^{k+1}$. Now, since $f_{\ov{c}}\neq 0$, the limit $\lim_{d\to \infty} (w^{d}_{(k+1)l})_{(k+1)\ov{c}}$ exists by assumption, and the sums $\sum_{n=(s+1)l}^d\binom{n}{s}g_nb_0^{n-s}$, $s=0,1,\ldots,k,$ converge, as $d\to \infty$, we have
$$\lim_{d\to \infty}\sum_{n=(k+1)l}^d\binom{n}{k+1}g_nb_0^{n-(k+1)}\in \K,$$ 
which proves the assertion.
\end{proof}

As we know from the one-variable case, the first coefficient of a fps under consideration plays an important role. In particular, it is important whether it is equal to zero or not. Therefore we introduce the following 
\begin{definition}
Let $f\in \X_q$, $q\in \bn$. We say that $f$ is a $q$-nonunit formal power series {if} $f_{(0,\ldots,0)}= 0$; otherwise we call it a q-unit.
\end{definition}

\begin{rmk}\textup{
Let us emphasize that in general the theory of formal power series is - roughly speaking - simpler if we restrict ourselves to $q$-nonunits. As an example one can mention here the composition of formal power series \textup{(}see e.g. \cite{GanM}, \cite{LoSt}\textup{)}. Namely, it is clear that necessary and sufficient conditions for the existence of composition of $g$ with $f$ are satisfied if $f$ is a q--nonunit fps. The same remark concerns the Right Distributive Law, which will be discussed in the next section. }
\end{rmk}

\section{General Right Distributive Law}\label{sec:4}

In this section we are going to consider formal power series of one variable only. Our goal is to remove ambiguities connected with the well known General Right Distributive Law. First, for convenience of the reader, we prove this law for arbitrary formal power series whose coefficients belong to $\R_+$.

\begin{theorem}\label{thm:23}
Let $A(x) = \sum\limits_{n=0}^\infty a_n x^n$, $B(x) = \sum\limits_{n=0}^\infty b_n x^n$ and $P(x) = \sum\limits_{n=0}^\infty p_n x^n$ be formal power series over $\R$ such that their coefficients belong to $\R_+$. The Right Distributive Law, that is,
\begin{equation}\label{eq:2.12}
( A \circ P ) ( B \circ P ) = ( A B ) \circ P
\end{equation}
holds if both $A \circ P$ and $B \circ P$ exist.
\end{theorem}

\begin{proof}
Suppose that both $A \circ P$ and $B \circ P$ exist. Then
\[
(A \circ P)(x) (B \circ P)(x) = \sum_{m=0}^\infty r_m x^m
\]
is a formal power series in $\X(\R)$ and $r_m \in \R_+$, $m\in \N_0$. Let
\[
P^n (x) = \sum_{k=0}^\infty p_k^{(n)} x^k,  \quad\text{$ n \in \N_0 $.}
\]
Since $P^{n+i} = P^nP^i$, the formula of the Cauchy product provides
\[
p_m^{(n+i)} = \sum_{j=0}^m \, p_j^{(n)} p_{m-j}^{(i)}
\]
for each $m \in \N_0.$ Then
\begin{align*}
 r_m & = \sum_{j=0}^m\bigg ( \sum_{n=0}^\infty a_n p_j^{(n)}\bigg ) 
\bigg( \sum_{i=0}^\infty b_i p_{m-j}^{(i)}\bigg ) \\
& = \sum_{n=0}^\infty a_n \left( \sum_{i=0}^\infty b_i \, \bigg(\sum_{j=0}^m p_j^{(n)}
\, p_{m-j}^{(i)} \bigg)\right)
 =  \sum_{n=0}^\infty a_n \left(\sum_{i=0}^\infty b_i \, p_m^{(n+i)}\right).
\end{align*}
Hence, by the absolute convergence of the series under consideration,
\begin{align*}
r_m=\sum_{j=0}^\infty \, a_j \left(\sum_{n=0}^\infty \, b_n p_m^{(n+j)}\right)
 = \sum_{j=0}^\infty a_j \left( \sum_{n=j}^\infty \, b_{n-j} p_m^{(n)}\right)
 =  \sum_{n=0}^\infty \,\bigg ( \sum_{j=0}^n a_j b_{n-j} \bigg) p_m^{(n)}
\end{align*}
for any $m\in \N_0$. However, the last term on the right--hand side is exactly the $m$-th coefficient of $( A B ) \circ  P.$ Thus, the composition $( A B ) \circ  P$ exists and equation (\ref{eq:2.12}) is satisfied.
\end{proof}

\begin{rmk}\textup{
Let us emphasize that in the above theorem we have considered fps whose coefficients belong to $\R_+$. In \cite{Gan1} the Author claims that the Generalized Right Distributive Law holds for fps over $\mathbb R$ while in \cite{GanM2} and \cite{GanM} he claims that it holds for fps over $\mathbb C$. The following example will show that unfortunately it is not a true claim.}
\end{rmk}

\begin{xmpl}\label {counterex}\textup{
Let $A=(a_0,a_1,a_2,\ldots)$ be a formal power series, where $a_n=\frac{(-1)^n}{\sqrt{n+1}},\, n\in \N_0$. The series $\sum_{n=0}^\infty a_n$ is conditionally convergent while the Cauchy product of the series with itself diverges (see \cite{Rud}, Example 3.49). Let $P=(1,0,0,\ldots)$. It is not difficult to see that $A\circ P=(\sum_{n=0}^\infty a_n,0,0,\ldots)$, so the composition $A\circ P$ exists and thus the (Cauchy) product of formal power series $(A\circ P)(A\circ P)$ is well-defined. But $(AA)\circ P$ does not exist, since the first coordinate of it, $((AA)\circ P)_0=\sum_{n=0}^\infty c_n$, where $c_n=a_0a_n+a_1a_{n-1}+\ldots+a_na_0$, diverges. This example invalidates the generality of the Right Distributive Law. }
\end{xmpl}

It appears that, adding the assumption of the existence of the composition $AB\circ P$ in Theorem \ref{thm:23}, one can prove the following version of the General Right Distributive Law.

\begin{thm}\label{thm:rdl1}
Let $A(z)=\sum_{n=0}^\infty a_nz^n,\, B(z)=\sum_{n=0}^\infty b_nz^n,\, P(z)=\sum_{n=0}^\infty p_nz^n$ be formal power series over $\mathbb{C}$. If $A\circ P$, $B\circ P$, and $(AB)\circ P$ exist, then the Right Distributive Law holds for $A,B$ and $P$, that is,
$$(A\circ P)(B\circ P)=(AB)\circ P.$$
\end{thm}

\begin{proof}
Since $A\circ P$, $B\circ P$ and $(AB)\circ P$ exist, for $m\in \N_0$, the following series, that are $m$-th coordinates of $A\circ P$, $B\circ P$ and $(AB)\circ P$, respectively,
$$(A\circ P)_m=\sum_{n=0}^\infty a_np^{(n)}_m,\,\quad (B\circ P)_m=\sum_{n=0}^\infty b_np^{(n)}_m,\,\quad ((AB)\circ P)_m=\sum_{n=0}^
\infty c_np^{(n)}_m,$$
where $p^{(n)}_m:=(P^n)_m $, $c_n=\sum_{j=0}^na_jb_{n-j}$, are convergent in $\C$.

For any $x\in(0,1]$ and $m\in \N_0$, we denote $$(A\circ P)_m(x):=\sum_{n=0}^\infty a_np^{(n)}_mx^n,\,\quad (B\circ P)_m(x):=\sum_{n=0}^\infty b_np^{(n)}_mx^n,\,\quad((AB)\circ P)_m(x):=\sum_{n=0}^\infty c_np^{(n)}_mx^n.$$
Observe that the series $(A\circ P)_m(x),\,(B\circ P)_m(x),\,((AB)\circ P)_m(x)$ are absolutely convergent power series for $x\in (0,1)$, since $(A\circ P)_m(1)=(A\circ P)_m,\,(B\circ P)_m(1)=(B\circ P)_m,\,((AB)\circ P)_m(1)=((AB)\circ P)_m$. By the Mertens Theorem (\cite{Rud}, Th. 3.50), it follows that for $s,t\in \N_0$, $x\in (0,1)$, the Cauchy product of series $(A\circ P)_s(x)$ and $(B\circ P)_t(x)$ is a convergent series $R^{s,t}(x)=\sum_{n=0}^\infty r^{s,t}_nx^n$, where $r^{s,t}_n:=\sum_{k=0}^na_kb_{n-k}p^{(k)}_sp^{(n-k)}_t$ and $(A\circ P)_s(x) (B\circ P)_t(x)=R^{s,t}(x)$.
Hence, for $m\in \N_0$, $x\in (0,1)$, we have
\begin{multline*}\sum_{j=0}^m(A\circ P)_{j}(x) (B\circ P)_{m-j}(x)=\sum_{j=0}^m R^{j,m-j}(x)=\sum_{j=0}^m\sum_{n=0}^\infty r^{j,m-j}_nx^n=\\
\sum_{j=0}^m \sum_{n=0}^\infty (\underbrace{\sum_{k=0}^na_kb_{n-k}p^{(k)}_jp^{(n-k)}_{m-j}}_{=r^{j,m-j}_n})x^n=\sum_{n=0}^\infty \sum_{j=0}^m (\sum_{k=0}^na_kb_{n-k}p^{(k)}_jp^{(n-k)}_{m-j})x^n=(\star),
 \end{multline*}
where the change of summation order is due to the convergence of series $\sum_{n=0}^\infty r^{j,m-j}_nx^n$.
Further,
\begin{multline*}(\star)=\sum_{n=0}^\infty\sum_{k=0}^n\sum_{j=0}^m (a_kb_{n-k}p^{(k)}_jp^{(n-k)}_{m-j})x^n=\sum_{n=0}^\infty\sum_{k=0}^n(\underbrace{\sum_{j=0}^m p^{(k)}_jp^{(n-k)}_{m-j}}_{p^{(n)}_m})a_kb_{n-k}x^n=\sum_{n=0}^\infty\sum_{k=0}^np^{(n)}_ma_kb_{n-k}x^n=\\\sum_{n=0}^\infty\underbrace{(\sum_{k=0}^na_kb_{n-k})p^{(n)}_m}_{c_n}x^n=((AB)\circ P)_m(x).
 \end{multline*}
Now, by Abel's Theorem (\cite{Rud}, Th. 8.2), we conclude that $\sum_{j=0}^m(A\circ P)_{j}(1) (B\circ P)_{m-j}(1)=((AB)\circ P)_m(1)$, which implies $((A\circ P)(B\circ P))_m=\sum_{j=0}^m(A\circ P)_{j}(B\circ P)_{m-j}=((AB)\circ P)_m$. Therefore, due to arbitrariness of $m\in \N_0$, the claim follows.
\end{proof}

\begin{remark}\textup{
In connection with the above result, we have to add one more explanation. In \cite{GanM2} {one can find} the Right Distributive Law (Theorem 5.6.2), where it is assumed that $A\circ P$ and $B\circ P$ exist. It is a consequence of Lemma 5.6.1 therein, which states that this assumption implies the existence of the composition $(AB)\circ P$. However, one has to notice that the proof of Lemma 5.6.1 is based on Theorem~5.4.6, in which there is the assumption that $deg(f)\neq 0$ (the degree of $f$ is nonzero). Unfortunately this assumption is missed both in Lemma 5.6.1 and Theorem 5.6.2. Finally, let us emphasize that in Example \ref{counterex} the degree of $A$ (which plays the same role as $f$ mentioned above) is equal to zero. It should explain all the ambiguities connected with this law. }
\end{remark}

\begin{rmk}\textup{
Theorem \ref{thm:rdl1} easily implies Abel's theorem asserting that if two series converge and their Cauchy product converges, then the product of their values equals the value of their Cauchy product.
To see this, consider series $\sum_{n=0}^\infty a_n,\,\sum_{n=0}^\infty b_n,\,\sum_{n=0}^\infty c_n$ such that $\sum_{n=0}^\infty c_n$ is the Cauchy product of $\sum_{n=0}^\infty a_n$ and $\sum_{n=0}^\infty b_n$ and the three series are convergent in $\C$. Let $P(x):=I(x)\in \X_1,\, A(x):=\sum_{n=0}^\infty a_nx^n,\, B(x):=\sum_{n=0}^\infty b_nx^n,\,C(x):=\sum_{n=0}^\infty c_nx^n$. Then $(A\circ P)=(\sum_{n=0}^\infty a_n)I$, $(B\circ P)=(\sum_{n=0}^\infty b_n)I$, $(AB)\circ P=(\sum_{n=0}^\infty c_n)I$, which by Theorem \ref{thm:rdl1} results in $(\sum_{n=0}^\infty a_n)(\sum_{n=0}^\infty b_n)=\sum_{n=0}^\infty c_n$.}
\end{rmk}

\section{General Multivariable Right Distributive Law}

In this section we will present analogous results to those of Section \ref{sec:4}, but now we are concerned with multivariable formal power series.

\begin{theorem}
	Let $A(y) = \sum_{n=0}^{\infty} a_n y^n$, $B(y)=\sum_{n=0}^{\infty}b_n y^n \in \X_1$ and $P(x) = \sum_{c\in C}P_cX^c \in \X_q$ be formal power series over $\R$ such that their coefficients belong to $\R_+$. If both $A\circ P$ and $B \circ P$ exist, then the General Multivariable Right Distributive Law holds, that is:
	\[ (A\circ P)(B\circ P) = (AB)\circ P. \]	
\end{theorem}
\begin{proof}% and $((AB)\circ P)(x) = \sum_{c\in C} g_cX^c$. We only need to prove that for every $c\in C$, $h_c = g_c$
	Let $(A\circ P)(B\circ P)(x) = \sum_{c\in C} h_cX^c$. By the existence of $A\circ P$ and $B\circ P$ and Definition \ref{df:composition}, we obtain
	\[ (A\circ P)(x) = \sum_{c\in C}(A\circ P)_cX^c = \sum_{c\in C}\left(\sum_{n=0}^{\infty}a_n p^{(n)}_c\right)X^c, \]
	\[ 
	(B\circ P)(x) = \sum_{c\in C}(B\circ P)_cX^c = \sum_{c\in C}\left(\sum_{n=0}^{\infty}b_n p^{(n)}_c\right)X^c,
	\]
	where $p^{(n)}_c:=(P^n)_c$.
	By definition of the (Cauchy) product, we get 
	\[ h_c = \sum_{l,m\in C:\, l+m=c} (A\circ P)_l (B\circ P)_m = \sum_{l,m\in C:\, l+m=c} \left( \sum_{j=0}^{\infty} a_j p^{(j)}_l\right) \left( \sum_{i=0}^{\infty} b_i p^{(i)}_m\right)=(\star), 
	\]
	and since the external sum has a finite number of components we can write 
	\[ (\star)= \sum_{j=0}^{\infty}a_j \left(\sum_{i=0}^{\infty}b_i \left( \sum_{l,m\in C:\, l+m=c}p^{(j)}_l p^{(i)}_m\right)\right) = \sum_{j=0}^{\infty} a_j \left(\sum_{i=0}^{\infty} b_i p^{(i+j)}_c\right).\]
	Now, by absolute convergence of considered series,
	\[ h_c = \sum_{j=0}^{\infty} a_j \left(\sum_{i=0}^{\infty}b_i p^{(i+j)}_c\right)=\sum_{j=0}^{\infty} a_j \left(\sum_{i=j}^{\infty} b_{i-j} p^{(i)}_c\right)= \sum_{i=0}^{\infty} \left( \sum_{j=0}^{i}a_j b_{i-j}\right)p^{(i)}_c\] for all $c\in C$. From this, definitions of the Cauchy product $AB$ and composition $AB\circ P$ we conclude that $(AB\circ P)_c=h_c$, $c\in C$. Therefore the composition $AB\circ P$ exists and the asserted equality holds true. \end{proof}

As we did in the previous section, we are going to prove that by adding the {assumption of the} existence of $(AB)\circ P$ we can extend the above theorem to the case of series with coefficients in $\C$.

\begin{theorem}\label{fin}
	Let $A(y) = \sum_{n=0}^{\infty} a_n y^n$, $B(y)=\sum_{n=0}^{\infty}b_n y^n \in \X_1$ and $P(x) = \sum_{c\in C}P_cX^c \in \X_q$ be formal power series over $\C$. Suppose that $A\circ P$, $B \circ P$ and $(AB)\circ P$ exist. Then 
	\[ 
	(A\circ P)(B\circ P) = (AB)\circ P, 
	\]
	that is, the General Multivariable Right Distributive Law holds.
\end{theorem}
\begin{proof}
	Using Definition \ref{df:composition}, we can write
	\[ (A\circ P) (x) = \sum_{c\in C} (A\circ P)_cX^c = \sum_{c\in C} \left( \sum_{n=0}^{\infty} a_n p^{(n)}_c\right)X^c \]
	and
	\[ (B\circ P) (x) = \sum_{c\in C} (B\circ P)_cX^c = \sum_{c\in C} \left( \sum_{n=0}^{\infty} b_n p^{(n)}_c\right)X^c, \]
	where $p^{(n)}_c:=(P^n)_c$.
	Now, using the definition of the Cauchy product we get
	\[ ((AB)\circ P) (x) = \sum_{c\in C}((AB)\circ P)_cX^c = \sum_{c\in C} \left( \sum_{n=0}^{\infty} c_n p^{(n)}_c \right)X^c, \]
	where $c_n =\sum_{k=0}^{n} a_k b_{n-k}$.\newline
	By assumption, the three compositions above exist, so for every $ c\in C$, the respective $c$-th coefficients are convergent series, that is, $((AB)\circ P)_c$, $(A\circ P)_c$, $(B\circ P)_c \in \C$. For any $t\in (0,1]$, let us define
	\[ (A\circ P)_c(t) :=  \sum_{n=0}^{\infty} a_n p^{(n)}_c t^n,\quad (B\circ P)_c(t):= \sum_{n=0}^{\infty} b_n p^{(n)}_c t^n, \quad ((AB)\circ P)_c(t) = \sum_{n=0}^{\infty} c_n p^{(n)}_c t^n. \]
	First, observe that $(A\circ P)_c(1) = (A\circ P)_c\in \C$, so for $t\in (0,1)$, $(A\circ P)_c(t)$ is absolutely convergent for any fixed $c\in C$. Similar properties are shared by $(B\circ P)_c(t)$ and $((AB)\circ P)_c(t)$, so we can apply the Mertens Theorem to deduce that, for $r,s\in C$, $t\in (0,1)$, the Cauchy product $((A\circ P)_r(t))((B\circ P)_s(t)) = Q^{r,s}(t)$ is convergent, where $Q^{r,s}(t) = \sum_{n=0}^{\infty}q_n^{r,s} t^n$ and $q_n^{r,s} = \sum_{k=0}^j a_k b_{n-k} p^{(k)}_r p^{(j-k)}_s.$ 
	
We want to prove that for every $c\in C$, $((A\circ P)(B\circ P))_c = ((AB)\circ P)_c$. To this end fix $c\in C$ and observe that for every $y\in(0,1)$ we have
	\[ \sum_{r,s\in C:\, r+s=c} (A\circ P)_r(t) (B\circ P)_s(t) = \sum_{r,s\in C:\, r+s=c}Q^{r,s}(t) = \sum_{r,s\in C:\, r+s=c}\sum_{n=0}^{\infty} q_n^{r,s} t^n = \]
	\[ = \sum_{r,s\in C:\, r+s=c} \sum_{n=0}^{\infty} \sum_{m=0}^{n}(a_m b_{n-m} p^{(m)}_r p^{(n-m)}_s)t^n = \sum_{n=0}^{\infty}\left( \sum_{m=0}^{n} a_m b_{n-m}\right)\left( \sum_{r,s\in C:\, r+s=c}p^{(m)}_r p^{(n-m)}_s\right) t^n =(\star),\]
	where the last equality holds because of the convergence of $Q^{r,s}(t)$. Then 
	\[ 
	(\star)= \sum_{n=0}^{\infty} c_n p^{(n)}_c t^j = ((AB)\circ P)_c(t).  
	\]
	Now, we refer to Abel's Theorem to conclude that for all $c\in C$ \begin{multline*}((A\circ P)(B\circ P))_c=\sum_{r,s\in C:\, r+s=c} (A\circ P)_r(A\circ P)_s=\\\sum_{r,s\in C:\, r+s=c} (A\circ P)_r(1) (A\circ P)_s(1)=((AB)\circ P)_c(1)= ((AB)\circ P)_c,\end{multline*} so that $((A\circ P)(B\circ P))_c = ((AB)\circ P)_c,$ for $c\in C$, which proves our claim.

\end{proof}


\begin{bibdiv}
\begin{biblist}

\bib{BM}{article}{
title={Further remarks on formal power series},
author={Borkowski, M.},
author={Ma\'ckowiak, P.},
journal={Comment. Math. Univ. Carolin.},
volume={53\normalfont{(4)}},
pages={549--555},
year={2012},
}

\bib{BG}{article}{
title={A note on formal power series},
author={Bugajewski, D.},
author={Gan, X.--X.},
journal={Comment. Math. Univ. Carolin.},
volume={51\normalfont{(4)}},
pages={595--604},
year={2010},
}

\bib{Daw}{article}{
title={The inverse and the composition in the set of formal Laurent series},
author={Bugajewski, D.},
journal={arXiv:2202.13948v1},
pages={1--30},
year={2022},
}

\bib{CMP}{article}{
title={Decomposition and eigenvectors of Riordan matrices},
author={Cheon, G.--S.},
author={Cohen, M.M. },
author={Pantelidis, N.},
journal={Linear Alg. Appl.},
volume={642},
pages={118--138},
year={2022},
}

\bib{Gan1}{article}{
title={A Generalized Chain Rule for formal power series},
author={Gan, X.--X.},
journal={Comm. Math. Anal.},
volume={2\normalfont{(1)}},
pages={37--44},
year={2007},
}

\bib{GanM}{book}{
   author={Gan, X.--X.},
   title={Selected Topics of Formal Analysis},
   series={Lecture Notes in Nonlinear Analysis}
   volume={15}
   publisher={Juliusz Schauder Center for Nonlinear Studies, Nicolaus Copernicus University},
   address={Toru\'n},
   date={2017},
}

\bib{GanM2}{book}{
   author={Gan, X.--X.},
   title={Formal Analysis. An Introduction},
   series={De Gruyter Studies in Mathematics}
   volume={80}
   publisher={De Gruyter},
   address={Berlin},
   date={2021},
}

%\bib{Hen}{book}{
%   author={Henrici, P.},
%   title={Applied and Computational Complex Analysis, vol.1},
%   series={Wiley Classics Library}
%   volume={43}
%   publisher={John Wiley and Sons},
%   date={1988},
%}

\bib{GK}{article}{
title={On composition of formal power series},
author={Gan, X.--X.}, 
author={Knox, N.},
journal={Int. J. Math. Math. Sci.},
volume={30\normalfont{(12)}},
pages={761--770},
year={2002},
}

\bib{Hau}{article}{
title={Formal power series in several variables},
author={Haukkanen, P.},
journal={Notes Number Theory Discrete Math.},
volume={25\normalfont{(4)}},
pages={44--57},
year={2019},
}

\bib{JLN}{article}{
title={Some algebraic structure of the Riordan group},
author={Jean--Louis, C.},
author={Nkwanta, A.},
journal={Linear Alg. Appl.},
volume={438\normalfont{(5)}},
pages={218--235},
year={2013},
}

\bib{LoSt}{book}{
   author={{\L}ojasiewicz, S.}, 
	 author={Stasica, J.},
   title={Formal Analysis and Analytic Functions},
   publisher={Jagiellonian University Press},
   address={Cracow},
   date={2005},
	 language={in Polish},
}

\bib{PraSp}{article}{
title={Unique summing of formal power series solutions to advanced and delayed differential equations},
author={Parvica, D.},
author={Spurr, M.},
journal={Discr. Cont. Dyn. Syst.,  Supplement Volume},
pages={730--737},
year={2005},
}

\bib{Rud}{book}{
   author={Rudin, W.},
   title={Principles of Mathematical Analysis},
   series={International Series in Pure and Applied Mathematics},
   publisher={McGraw--Hill Inc.},
   date={1976},
}

\bib{Sib}{article}{
title={Formal power series solutions in a parameter},
author={Sibuya, Y.},
journal={J. Diff. Eq.},
volume={190\normalfont{(2)}},
pages={559--578},
year={2003},
}

\bib{S}{book}{
author={Stanley, R.P.},
   title={Enumerative Combinatorics},
   series={Cambridge Studies in Advanced Mathematics},
	 volume={49},
   publisher={Cambridge University Press},
   date={1997},
}

\end{biblist}
\end{bibdiv}
\end{document}